\documentclass[11pt,a4paper]{article}
\usepackage{geometry}
\usepackage[british]{babel} 
\usepackage[T1]{fontenc}
\usepackage[utf8]{inputenc}
\usepackage{color}
\usepackage{url}
\usepackage{fourier,charter}
\usepackage{varioref}\labelformat{equation}{(#1)}

\usepackage[round]{natbib}
\usepackage{enumerate,mathrsfs}
\usepackage[modulo,mathlines]{lineno}
\newcommand{\cov}[2][]{
   \ensuremath{\textrm{cov}_{#1}\!\left\{\displaystyle{#2}\right\}}}
\renewcommand{\d}{\,\mathrm{d}}

\newcommand{\equald}{\stackrel{d}{=}}
\newcommand{\E}[2][]{
   \ensuremath{\mathbb{E}_{#1}\!\left\{\displaystyle{#2}\right\}}}

\newcommand{\inv}{^{-1}}
\newcommand{\pr}[2][]{
   \ensuremath{\mathbb{P}_{#1}\!\left\{\displaystyle{#2}\right\}}}
\newcommand{\pd}[2]{\displaystyle\frac{\partial #1}{\partial #2}}   
\newcommand*{\QED}[1][\textsc{\small qed}]{%
 \leavevmode\unskip\penalty9999 \hbox{}\nobreak\hfill\quad\hbox{#1}}   
\newcommand{\Real}{\mathbb{R}}
\newcommand{\tr}{\mbox{\rm tr}}
\newcommand{\T}{^{\top}}
\newcommand{\var}[2][]{
   \ensuremath{\textrm{var}_{#1}\!\left\{\displaystyle{#2}\right\}}}
\renewcommand{\vec}{\mathop{\mathrm{vec}}\nolimits}
\newtheorem{theorem}{Theorem}
\newtheorem{lemma}[theorem]{Lemma}

\newtheorem{proposition}[theorem]{Proposition}
\newcommand{\N}{\mathrm{N}{}}

\newcommand{\SUN}{\mathrm{SUN}{}}
\let\phi=\varphi
\title{\bf Some properties of the unified skew-normal distribution}
\author{
   {\bf Reinaldo B. Arellano-Valle}\\
    Departamento de Estadística\\
    Pontificia Universidad Católica de Chile\\
    Chile
  \and
   {\bf  Adelchi Azzalini} \\
    Dipartimento di Scienze Statistiche \\
    Università di Padova\\
    Italia
 }
\date{\today}
\begin{document}

\maketitle
\begin{abstract}
For the family of multivariate probability distributions variously denoted
as unified skew-normal, closed skew-normal and other names, a number of
properties are already known, but many others are not, even some basic
ones. The present contribution aims at filling some of the missing gaps.
Specifically, the moments up to the fourth order are obtained,
and from here the expressions of the Mardia's measures of multivariate
skewness and kurtosis. Other results concern the property of log-concavity
of the distribution, and closure with respect to conditioning on intervals.
\end{abstract}
\vspace{2ex}

\section{The unified skew-normal distribution}
\subsection{Early development, applications and some open problems}  \label{s:intro}

In recent years, there has been a vigorous impulse in the development
of flexible parametric families of distributions.
This activity is specially lively and stimulating in the multivariate setting,
correspondingly to the  ever increasing availability and
treatment of multivariate data in applied work.

An active direction of research within this process is represented by a family of
continuous distributions which has originated as a generalization of the
multivariate skew-normal (SN) distribution, which itself is a generalization
of the classical normal distribution; for a review of the SN
distribution and its ramifications, see \cite{azza:capi:2014}.
The generalization we are concerned with has originated from multiple
independent sources, with some differences in the technical development,
but with common underlying structure, as explained later.
Specifically, \cite{gonz:domi:guptAK:2004} and \cite{gonz:etal:2004-inSE}
have developed the `closed skew-normal distribution'.
Motivated by Bayesian inference considerations,
\cite{lise:lope:2003} have presented the `hierarchical skew-normal'.
Another related construction is the `fundamental skew-normal'
proposed by \cite{arel:gent:2005}, who also consider a second version
of closed skew-normal.

The interconnections among these apparently  separate formulations
have been examined by  \cite{arel:azza:2006}, showing their essential
equivalence, as well as the presence of overparameterizations
in some cases.
To accomplish their project, they introduced a unifying version which
embraces the above-recalled specific proposals,
removing at the same time the existing overparameterizations.
This version was hence denoted `unified skew-normal (SUN)  distribution'.
Its main formal properties  will be summarized in the next subsection.
However, in essence, the constructive mechanism  starts from
a $(d+m)$-dimensional normal distribution, where $m$ of the components
play a role of hidden variables which modify
non-linearly the remaining $d$ components via the presence of
a certain conditioning event on the hidden components.
The construction leads to a $d$-dimensional non-normal distribution,
with the regular normal distribution included as a special case.
We shall refer to this  distribution as a $\SUN_{d,m}$.

The SUN family constitutes a superset of the SN family,
more specifically  the so-called `extended skew-normal (ESN)
family', to which the SUN family reduces if $m=1$.
Its building mechanism based on $m$ latent variables leads to
certain properties not amenable to the SN and ESN distribution.
An important specific fact is closure of the family
with respect to convolution; specifically, the sum of two independent
SUN variables of type $\SUN_{d,m_1}$ and $\SUN_{d,m_2}$ is of type
$\SUN_{d,m_1+m_2}$.
This property has proved convenient in a number of operational
formulations which employ the  SUN distribution as its core
stochastic component.

The closed skew-normal and the SUN distributions have been applied in
a wide range of applied domains, and their relevance appears to be growing.
The following is a non-exhaustive list of methodologies and applied
domains  where these distributions have been employed:
stochastic frontier analysis in the context of productivity analysis,
considered by \cite{domi:gonz:ramo:guptAK:2007}, \cite{colombi:2013},
\cite{colo:kumb:etal:2014}, \cite{kumb:laiHP:2016};
various models for the analysis of spatial data have been
introduced by \cite{alla:nave:2007}, \cite{hoss:eids:moha:2011},
\cite{kari:moha:2012},  \cite{rims:omre:2014}, among others;
analysis of longitudinal data for the distribution of random effects
in work of \cite{ghal:zadk:2019}, and again \cite{colombi:2013};
combination of phase~II and~III clinical trials,  by \cite{azza:bacc:2010};
seismic inversion methodology for geological problems,
by \cite{kari:omre:moha:2010} and \cite{reza:eids:muke:2014};
extended formulations of Kalman filter by \cite{kimHM:ryu:mall:gent:2014}
and \cite{reza:eids:2016};
application to small area estimation by  \cite{dial:raoJNR:2018}.
In the context of binary data,  \cite{durante:2019} has shown that,
under Gaussian priors for the probit coefficients, the posterior
distribution  has an exact unified skew-normal distribution;
this formulation lends itself to interesting developments,
such as those of \cite{fasa:dura:zane:2019} and \cite{fasa:dura:2020}.

While the SUN distribution is mathematically quite tractable and it
enjoys a number of  appealing formal properties, it is inevitably
more complex than its progenitor, that is, the skew-normal distribution.
Consequently there are several aspects which are still unexplored,
or only partly explored; this situation concerns even some rather 
basic properties.
A case in point is represented by the computation of the moments and
associated quantities,  of which little is known at present,
as we shall discuss in more detail later on.
This problem represents the main target of the present contribution,
tackled in Section~\ref{s:moments}.
Additional properties are examined in Section~\ref{s:other-properties},
namely the study of the log-concavity of the density and
the conditional distribution of a SUN variable when some
of its components belong to a given interval.

\subsection{Main properties of the SUN family}   \label{s:properties}

We summarize the main facts about the SUN family; this term is used
to embrace also the closed skew-normal and other essentially equivalent
classes, provided a suitable parameterization is adopted.
The notation here is the one of Subsection~7.1.2 of \cite{azza:capi:2014},
which is largely the same of \cite{arel:azza:2006}, with minor variations.

For positive integers $d$ and $m$, consider the $(d+m)$-dimensional
normal random variable
\begin{equation}
     \pmatrix{X_0\cr X_1} \sim \N_{d+m}\left(0, \Omega^* \right) \,,
         \qquad
        \Omega^* = \pmatrix{ \bar\Omega & \Delta \cr
                             \Delta\T &  \Gamma },
  \label{e:(d+m)-normal}
\end{equation}
where $\Omega^*$ is a full-rank correlation matrix.
Define $Z$ to be a $d$-dimensional random variable with the same
distribution of  $(X_0|X_1+\tau>0)$,
where $\tau=(\tau_1,\ldots,\tau_m)\T$ and the notation $X_1+\tau>0$
means that the inequality sign must hold component-wise for each one
of the $m$ components.
Next, introduce the transformed variable $Y=\xi+\omega\,Z$,
where $\xi=(\xi_1,\ldots,\xi_d)\T$ and
$\omega$ is a $d\times d$ diagonal matrix with positive diagonal elements
$\omega_1,\ldots, \omega_d$, and denote $\Omega=\omega\bar\Omega\omega$.
It can be show that the density of $Y$ at $x\in\Real^d$ is
\begin{equation}    \label{e:sun-pdf}
  f_Y(x)= \phi_d(x-\xi;\Omega)\:
     \frac{\Phi_m\left\{\tau +  \Delta\T\bar\Omega\inv\omega\inv (x-\xi);
          \Gamma - \Delta\T\bar\Omega\inv\Delta\right\}}{\Phi_m(\tau;\Gamma)}
\end{equation}
where  $\phi_h(u;\Sigma)$ and  $\Phi_h(u;\Sigma)$  denote the $\N_h(0,\Sigma)$
density function and distribution function at $u\in\Real^h$, respectively,
for any symmetric $(h \times h)$ positive-definite matrix $\Sigma$.
In this case, we shall write $Y \sim \SUN_{d,m}(\xi, \Omega, \Delta, \tau, \Gamma)$.

The SUN family enjoys numerous formal properties. For instance,
we have already anticipated in Subsection~\ref{s:intro} that this family is
closed with respect to convolution. Many other interesting facts hold,
but it would take too much space to review all such properties here,
and we only recall those which are required for the subsequent development;
additional information is summarized in Section~7.1 of \cite{azza:capi:2014}.
A key fact is the expression of the moment generating function, $M(t)$  or,
equivalently, the cumulant generating function of \ref{e:sun-pdf}
as given by \cite{arel:azza:2006} is
\begin{equation}    \label{e:sun-cumul}
  K(t) = \log M(t) = \xi\T t + 2\inv t\T\Omega t \:
         + \log\Phi_m(\tau+\Delta\T\omega t; \Gamma) - \log\Phi_m(\tau; \Gamma),
   \hspace{3em}t\in\Real^d;
\end{equation}
essentially as in \citet{gonz:domi:guptAK:2004} and \citet{gonz:etal:2004-inSE},
up to a change of parameterization.
From this expression, many other results can be derived. One of them
is represented by the rule for obtaining the distribution of an affine
transformation:
if $a$ is a $p$-vector and $A$ is a full-rank $d\times p$ matrix, then
\begin{equation}
  a+ A\T Y \sim \SUN_{p,m}(a+A\T\xi, A\T\Omega A, \Delta_A, \tau, \Gamma)
  \label{e:sun-affine}
\end{equation}
where $\Delta_A= \mathrm{Diag}(A\T \Omega A)^{-1/2} A\T\omega\Delta$,
using the notation $\mathrm{Diag}(M)$  to denote the diagonal matrix
formed by the diagonal elements of a square matrix $M$,
as in \citet[p.\,455]{mard:kent:bibb:1979}.
Clearly, \ref{e:sun-affine} can be used to compute the distribution of
$p$-dimensional marginals.

Another result to be used in our development is the expression of
the distribution function, which has been given in Lemma 2.2.1 of
\citet{gonz:etal:2004-inSE}.
Since we adopt the SUN formulation for the reasons discusses by \cite{arel:azza:2006},
we shall use the equivalent expression, given by \cite{azza:bacc:2010},
\begin{equation}
\label{e:sun-cdf}
  F_Y(y)= \pr{Y\le y} =\frac{\Phi_{d+m}(\tilde{z}; \tilde\Omega)}{\Phi_m(\tau; \Gamma)}
\end{equation}
where
\[
    \tilde{z} = \pmatrix{\omega\inv(y-\xi) \cr \tau }, \qquad
    \tilde\Omega= \pmatrix{\bar\Omega & -\Delta \cr -\Delta\T & \Gamma}\,.
\]

There exist two stochastic representations of the SUN distribution, or equivalently
two constructive ways to generate a random variable $Y$ with density \ref{e:sun-pdf}.
The first of these is essentially the above-described process leading from
the normal variable $X$ in \ref{e:(d+m)-normal} to the variable $Y$, via the
intermediate variable $Z$. This is denoted `representation by conditioning'
since it operates through the condition  $X_1+\tau>0$.

The other stochastic representation is of convolution type, that is,
as the distribution of the sum of two independent random variables.
Specifically, from the above-defined quantities,
introduce $\bar\Psi_\Delta =\bar\Omega-\Delta\Gamma\inv\Delta\T$,
and the two  independent variables  $U_0\sim\N_d(0,\bar\Psi_\Delta)$
and $U_{1,-\tau}$  which is obtained  by the component-wise truncation
below $-\tau$ of a variate  $U_1\sim\N_m(0,\Gamma)$.
Then, $Y \sim \SUN_{d,m}(\xi, \Omega, \Delta, \tau, \Gamma)$
can be  expressed  via the so-called additive representation
\begin{equation}    \label{e:sun-additive}
    Y \equald \xi+\omega\:\left(U_0 + \Delta\Gamma\inv\,U_{1,-\tau}\right)
\end{equation}
which will play a key role in our development.
For a detailed discussion of the interplay of these two stochastic
representations, see Section~2.1 of \cite{arel:azza:2006}.

Although the moment generating function $M(t)$ has been known since the
early work on this theme,  it  has not translated into decisive
advances in the computation of moments and  cumulants.
Most of the available results for $\E{Y}$ and $\var{Y}$ are limited
is some way or another. For instance, results in Section~3 of
\cite{guptAK:gonz:domi:2004} refer to the case $m=d$, and even so
they employ very involved auxiliary functions.
For the case where $\Gamma$ is a diagonal matrix,
\cite{arel:azza:2006} provide explicit expressions  for the expected
value and the variance matrix, applicable for all $d$ and $m$.

To our knowledge, the general expression of $\E{Y}$ has been obtained
by \cite{azza:bacc:2010}.
This expression involves the following quantities:
$\tau_{-j}$ denotes the vector obtained by removing the $j$ component of $\tau$,
  for $j=1,\dots,m$;
$\Gamma_{-j}$ is the $(m-1)\times(m-1)$
  matrix obtained by removing the $j$th row and column of $\Gamma$;
$\gamma_{-j}$ denotes the $j$th column of $\Gamma_{-j}$;
finally, $\tilde\Gamma_{-j}= \Gamma_{-j} - \gamma_{-j}\gamma_{-j}\T$.
Then the mean value can be written as
\begin{equation}
  \E{Y} = \left.\frac{\d K(t)}{\d t}\right|_{t=0}
        = \xi + \omega\,\Delta\, \frac{1}{\Phi_m(\tau; \Gamma)}\:\nabla \Phi_m
  \label{e:sun-mean}
\end{equation}
where $\nabla \Phi_m$ is the $m$-vector with $j$th element
\begin{equation}
 (\nabla \Phi_m)_j = \cases{\phi(\tau_j) & if $m=1$, \cr
    \phi(\tau_j)\:\Phi_{m-1}\left(\tau_{-j}-\Gamma_{-j}\tau_j;\tilde\Gamma_{-j}\right) &
     if $m>1$.}
  \label{e:Phi-derivatives}
\end{equation}

An expression of type \ref{e:sun-mean} or \ref{e:Phi-derivatives}
can be regarded as `essentially explicit', at least for moderate
values of $m$, even if it involves the  distribution function
of a multivariate normal distribution function, $\Phi_m$.
The phrase `essentially explicit'
seems justified in  the light of the current advances for computing
$\Phi_m$,   similarly to the process which, a few decades ago,
has led to consider `explicit' an expression  involving the
univariate normal distribution function, $\Phi$.

Some intermediate expressions of the SUN variance matrix have been provided
by \cite{guptAK:aziz:2012} and \cite{guptAK:aziz:ningW:2013}, where the word
`intermediate' reflects the presence in their result of the matrix
of the second derivatives of $\Phi_m$.
Since these second derivatives have been provided in an explicit form only
for some special sub-cases of the SUN family, the question of the general
expression of the SUN variance matrix appears to be open.
This is the problem to be tackled in our next section, followed by
consideration of higher order moments.

\section{Moments and related quantities} \label{s:moments}
\subsection{The variance matrix}
We compute the variance matrix $\var{Y}$ using second-order differentiation
of the cumulant generating function \ref{e:sun-cumul}.
Write
\begin{equation}
   \frac{\d K(t)}{\d t}
    = \xi + \Omega t + \frac{\d \log P(t)}{\d t}
    = \xi + \Omega t + \frac{1}{P(t)}\: \frac{\d P(t)}{\d t}
   \label{e:sun-cumul-grad}
\end{equation}
where $P(t) =  \Phi_m(\tau+ \Delta\T\omega t; \Gamma)$.
The only non-obvious terms are $\partial P/\partial t_j$, for $j=1,\dots,m$.
Denote
\[  u_j= (\tau+ \Delta\T\omega t)_j
        =\tau_j+\Delta_j\T\omega t,
\]
where  $\Delta_j$ is the  $j$th column of $\Delta$, for $j=1,\ldots,m$.
For notational simplicity, we focus on $j=1$ since the other terms are analogous.
Write the joint $m$-normal density involved by $P$ as the product of the
first marginal component times the conditional density of the other components,
leading to
\begin{equation}
  P(t) = \int_{-\infty}^{u_1}\cdots \int_{-\infty}^{u_m}    \phi(x_1) \:
       \phi_{m-1}(x_{-1} -\mu_{-1}(x_1); \tilde\Gamma_{-1})\: \d x_1 \d x_{-1}
  \label{e:P1}
\end{equation}
where  $x_{-1}$ is the $(d-1)$-vector obtained by removing the first
component of $x$,  $\mu_{-1}(x_1)= \gamma_{-1} x_1$ denotes the mean value
of the conditional normal  distribution  when the first component
of $\N_m(0,\Gamma)$ is fixed at $x_1$,  and
$\tilde\Gamma_{-1}$ denotes the corresponding variance matrix;
we have used the quantities introduced in connection with \ref{e:sun-mean}.
Therefore
\begin{eqnarray}
 \frac{\partial P}{\partial t_1}
   =   \frac{\partial u_1}{\partial t_1}\:\frac{\partial P}{\partial u_1}
  &=& \left(\omega\Delta \right)_1 \: \phi(u_1)\:
      \int_{-\infty}^{u_2}\cdots \int_{-\infty}^{u_m} \:
     \phi_{m-1}(x_{-1} -\mu_{-1}(u_1); \tilde\Gamma_{-1})\:  \d x_{-1}
            \nonumber \\
  &=& \left(\omega\Delta \right)_1 \:  \phi(u_1)\:
     \Phi_{m-1}(\tau_{-1}- \gamma_{-1}u_1; \tilde\Gamma_{-1})
      \label{e:dP,t1}
\end{eqnarray}
where  the term $\Phi_{m-1}(\cdot{\color{red};}\cdot)$ must be interpreted as 1 when $m=1$.
This convention will apply also to subsequent expressions.

Application of \ref{e:dP,t1} with the other values of the subscript $j$
produces the entire gradient of $P$.
Next, evaluation of the gradient \ref{e:sun-cumul-grad} at $t=0$
delivers the mean vector \ref{e:sun-mean}.

The second derivative of $K(t)$ is obtained by differentiation of
\ref{e:sun-cumul-grad}, yielding
\begin{equation}
   \frac{\d^2 K(t)}{\d t\, \d t\T} = \Omega +
      \frac{\d \hspace{1em}}{\d t\T}\left( \frac{\d \log P(t)}{\d t}\right)
   \label{e:sun-cumul-hess}
\end{equation}
where two generic entries of the final term are of the type
\[
  \frac{\partial^2 \log P}{\partial t_1 \partial t_2} =
   - \frac{1}{P^2}\pd{P}{t_1}\pd{P}{t_2}  + \frac{1}{P} \pd{^2 P}{t_1\,\partial t_2}.
\]
The first summand on the right side is the product of quantities of type \ref{e:dP,t1}.
For the second summand  consider first the case with $t_1\not=t_2$,
and follow a similar logic used for \ref{e:P1},
but now separate out two components.
Focusing of the first two components, for notational simplicity, write
\begin{equation}
  P(t) =
    \int_{-\infty}^{u_1}\cdots \int_{-\infty}^{u_m}  \phi_2(x_{1:2}; \Gamma_{1:2}) \:
    \phi_{m-2}\left(x_{-(1:2)} -\mu_{-(1:2)}(x_{1:2}); \tilde\Gamma_{-(1:2)}\right)\:
    \d x_{1:2} \d x_{-(1:2)}
  \label{e:P2}
\end{equation}
where $x_{1:2}=(x_1,x_2)$, $\Gamma_{1:2}$ is the submatrix of $\Gamma$
formed by its  top-left $2\times 2$ block, and so on, in the same logic
and notational scheme used before.

Here we have implicitly assumed that $m\ge 2$.
This is a legitimate assumption since the case with $m=1$
corresponds to the ESN distribution, for which $\var{Y}$ has been given
by \cite{capi:azza:stan:2003} along with other moment-related results
of the ESN distribution.

The mixed derivative at $t_1=t_2=0$  is
\begin{eqnarray}
 \left.\frac{\partial^2 P(t)}{\partial t_1 \partial t_2}\right|_{t_1=0, t_2=0}
   &=& \left(\omega\Delta\right)_{1:2}\: \phi_2(\tau_{1:2}) \:
          \Phi_{m-2}(\tau_{-(1:2)} - \mu_{-(1:2)}(\tau_{1:2}); \tilde\Gamma_{-(1:2)})
       \left(\Delta\T \omega\right)_{1:2}
     \label{e:d2P,t1:t2}
\end{eqnarray}
where $\mu_{-(1:2)}(\tau_{1:2})$ denotes the conditional mean of the components $(3,\dots,m)$
conditionally on $x_{1:2}=\tau_{1:2}$ and $\tilde\Gamma_{-(1:2)}$ denotes the conditional
variance.
It must be intended that the term $\Phi_{m-2}(\cdot)$  is  1 when $m=2$.
Expression \ref{e:d2P,t1:t2} is immediately adapted to any two other components
$(t_j, t_k)$,  provided $j\not=k$.

When $j=k$, take $j=k=1$ for simplicity of notation and write
\[
  \frac{\partial^2 \log P}{\partial t_1^2}
   = \frac{\partial}{\partial t_1} \left(\frac{1}{P}\: \frac{\partial P}{\partial t_1}\right)
   = -\frac{1}{P^2}\left(\frac{\partial P}{\partial t_1}\right)^2
     + \frac{1}{P} \left(\frac{\partial^2 P}{\partial t_1^2}\right)
\]
where $(\partial P/\partial t_1)$ is given by \ref{e:dP,t1}.
Consider its core part $(\partial P/\partial u_1)$  and take  the successive derivative
\begin{eqnarray*}
 \frac{\partial^2 P}{\partial u_1^2}
     &=& \frac{\partial}{\partial u_1} \left(  \frac{d}{\partial u_1}
          \phi(u_1) \:\Phi_{m-1}(\tau_{-1}- \mu_1(u_1); \tilde\Gamma_{-1})\right) \\
     &=&  \frac{\partial}{\partial u_1} \:   \left(
          \phi(u_1) \:\Phi_{m-1}(\tau_{-1}- \gamma_{-1}u_1); \tilde\Gamma_{-1})\right) \\
     &=&
         - u_1 \phi(u_1)  \Phi_{m-1}(\tau_{-1}- \gamma_{-1}u_1); \tilde\Gamma_{-1})
        + \phi(u_1) \phi_{m-1}(\tau_{-1}- \gamma_{-1}u_1); \tilde\Gamma_{-1})
                1_{m-1}\T \frac{\partial (\tau_{-1}- \gamma_{-1}u_1) }{\partial u_1}
          \\
      &=&   - u_1 \phi(u_1)  \Phi_{m-1}(\tau_{-1}- \gamma_{-1}u_1); \tilde\Gamma_{-1})
          + \phi(u_1) \phi_{m-1}(\tau_{-1}- \gamma_{-1}u_1); \tilde\Gamma_{-1})
                1_{m-1}\T (- \gamma_{-1})   .
\end{eqnarray*}
Hence the second derivative $(\partial^2 P/\partial t_1^2)$ evaluated at $t_1=0$  is
\begin{eqnarray}
  \left.\frac{\partial^2 P}{\partial t_1^2}\right|_{t_1=0}
  &=&  (\omega\Delta)_{11} \:\big\{
        -\tau_1 \phi(\tau_1)\:\Phi_{m-1}(\tau_{-1}-\gamma_{-1}\tau_1;\tilde\Gamma_{-1})
         \nonumber  \\[-2ex]
   && \qquad\quad
     - \phi(\tau_1) \phi_{m-1}(\tau_{-1}- \gamma_{-1}\tau_1); \tilde\Gamma_{-1})
          1_{m-1}\T \gamma_{-1})    \big\}\:(\Delta\T\omega)_{11} \nonumber\\
  &=& (\omega\Delta)_{11} \:\big\{- \phi(\tau_1) \big[
       \tau_1\:\Phi_{m-1}(\tau_{-1}-\gamma_{-1}\tau_1;\tilde\Gamma_{-1}) \nonumber \\
  && \qquad\quad    + \phi_{m-1}(\tau_{-1}- \gamma_{-1}\tau_1; \tilde\Gamma_{-1}) \:
          1_{m-1}\T\gamma_{-1}\big]\big\}\: (\Delta\T\omega)_{11}
   \label{e:d2P,t1}
\end{eqnarray}
which, similarly to earlier expressions, must be replicated for the other values of $j$.

Finally, as a general expression encompassing all terms in a matrix notation,
we arrive at
\begin{equation}
  \var{Y} = \Omega + \omega\Delta H \Delta\T \omega = \Sigma,
  \label{e:sun-var}
\end{equation}
say, where $H$ is the matrix formed by the elements other than $\Delta\omega$
given in \ref{e:dP,t1}, \ref{e:d2P,t1:t2}  and \ref{e:d2P,t1}.
Even if  we have derived \ref{e:sun-var} under the assumption that $m\ge2$,
a subsequent inspection has shown that the expression remains valid provided
the above derivatives are computed setting the $\phi_{m-1}$ and $\Phi_{m-1}$
terms equal to 1 when $m=1$, as we recover the known expressions
for the ESN distribution.
With this convention, \ref{e:sun-var} holds for all $m$.

Starting from a different motivation, expressions for the
derivatives of $\Phi_m$ similar to those obtained above have been presented
in Lemma~2.3 of \cite{arel:cast:losc:2013}.
Their motivation was the computation of mean value and the variance matrix
of the  truncated multivariate normal distribution,  which are given in
their Lemma~2.2.
Taking into account the additive representation \ref{e:sun-additive} of a
SUN variable, those expression could also be used to derive the SUN
lower moments.


\subsection{Higher-order moments}

While in principle one could consider successive differentiations of $K(t)$
to compute higher-order moments, this process becomes algebraically cumbersome.
We therefore follow another route, based on the additive representation
\ref{e:sun-additive}.

Our plan of work is as follows.
A preliminary step is the development of various expressions concerning
moments of the sum of two independent random vectors, which are
presented separately in an appendix.
Simplification can be obtained by the using the fact that one of the
components of \ref{e:sun-additive} is a zero-mean normal variable.
On this front, we benefit from  the  extensive literature on
computational method for the moments of a multivariate truncated normal distribution.
Combining  representation \ref{e:sun-additive} with these results for the moments
of a truncated normal variable, we obtain expressions for the desired
SUN moments.

For a $p$-dimensional random variable $X$, define its moments up to the fourth
order as
\begin{eqnarray*}
\mu_1(X)&=&\E{X},\\
\mu_2(X)&=&\E{X\otimes X\T}=\E{XX\T},\\
\mu_3(X)&=&\E{X\otimes X\T\otimes X}=\E{X\otimes XX\T}=\E{XX\T\otimes X}=\E{\vec(XX\T)X\T},\\
\mu_4(X)&=&\E{X\otimes X\T\otimes X\otimes X\T}=\E{XX\T\otimes XX\T}=\E{\vec(XX\T)\vec(XX\T)\T}.
\end{eqnarray*}
provided the involved expected values exist.
The equivalence of the various expressions for a given moment follows
from standard properties of the Kronecker product.
The $\vec$ operator stacks the columns of a matrix in a single vector.

Also, the following notation will be used, adopted from
\cite{magn:neud:1979as} and \cite{neud:wans:1983}.
For arbitrary natural numbers $s, p, q$,
denote by  $e_{i:s}$  the  $i$th $s$-dimensional unit vector
formed by all 0's except a 1 in the $i$th position,
and from here define $E_{ij}=e_{i:p}e_{j:q}\T$.
Further, denote by
\[
    K_{pq}=\sum_{i=1}^p\sum_{j=1}^qE_{ij}\otimes E_{ij}\T
\]
the $pq$-dimensional square commutation matrix and let $K_r=K_{rr}$.

For algebraic convenience, we rewrite \ref{e:sun-additive} in an equivalent form.
Introduce the quantities
\begin{equation}
 \Lambda=  \omega\Delta\Gamma\inv, \qquad
   \Psi = \Omega- \omega \Delta\Gamma\inv\Delta\T \omega
  \label{e:Lambda,Psi}
\end{equation}
and denote by $\Psi^{1/2}$ the unique symmetric square root of $\Psi$;
however, it would make no difference if another square root of $\Psi$
is considered.
The fact that $\Psi>0$ follows from the assumption that $\Omega^*$ in
\ref{e:(d+m)-normal} has full rank.
Given a $(d+m)$-dimensional variable
\begin{equation}
  Z_0 = \pmatrix{V \cr W} \sim \N_{d+m}\pmatrix{I_d& 0 \cr 0 & \Gamma }
  \label{e:(d+m)-normal'}
\end{equation}
denote $U\equald (W\mid W+\tau>0)$, so that \ref{e:sun-additive} becomes
\begin{equation}
   Y \equald  \xi+ X = \xi + \Lambda\,U + \Psi^{1/2}\,V \,.
   \label{e:sun-additive2}
\end{equation}

\begin{proposition}[SUN moments]   \label{th:sun-moments}
Consider $X=\Lambda U+ \Psi^{1/2}V\sim \SUN_{d,m}(0,\Omega,\Delta,\tau,\Gamma)$,
where $U$, $V$ and other involved quantities are defined in connection with
expressions \ref{e:Lambda,Psi}--\ref{e:sun-additive2}. Then:
\begin{eqnarray*}
\mu_1(X)&=&\Lambda \mu_1(U),\\
\mu_2(X)&=&\Lambda \mu_2(U)\Lambda \T  + \Psi,\\
\mu_3(X)&=&(\Lambda\otimes\Lambda)\mu_3(U)\Lambda\T + (I_{d^2}+K_d)(\Lambda \mu_1(U)\otimes \Psi)
  + \vec(\Psi)\mu_1(U)\T\Lambda\T,\\
\mu_4(X)&=&(\Lambda\otimes\Lambda)\mu_4(U)(\Lambda\otimes\Lambda)\T
   + \vec(\Psi)\vec(\Psi)\T \\
&& + (I_{d^2}+K_d)\big\{(\Lambda \mu_2(U)\Lambda\T\otimes \Psi)
   + (\Psi\otimes \Lambda \mu_2(U)\Lambda\T) +  \Psi\otimes \Psi\big\}\\
&& + \vec(\Lambda \mu_2(U)\Lambda\T)\vec(\Psi)\T
   + \vec(\Psi)\vec(\Lambda \mu_2(U)\Lambda\T)\T \,.
\end{eqnarray*}
Moreover, the variance matrix of  $X$ is
\begin{equation}
    \var{X}=\Omega-\Lambda(\Gamma-\Sigma_U)\Lambda\T = \Sigma,
    \label{e:sun-var2}
\end{equation}
say, having set $\Sigma_U=\var{U}=\mu_2(U)-\mu_1(U)\mu_1(U)\T$.
The inequality $\Omega > \Sigma$ holds, in the sense that the
difference of the matrices is positive definite.
\end{proposition}
\noindent\emph{Proof.}
Make use of Proposition~\ref{th:moments-of-sum0} in the Appendix for the moments
of a linear combination of two independent multivariate variables,
combined with expressions in Proposition~\ref{th:moments-normal} for the
moments of $V$ under the assumption  $V\sim\N_d(0, I_d)$.
After some algebraic simplifications, one arrives at the stated expressions.
The term $\Gamma-\Sigma_u$ of \ref{e:sun-var2}
is a positive  definite matrix, taking account \ref{e:var-by-selection-ineq}
in the appendix, which in the present case holds in the strict version of
the matrix inequality. This implies that $\Omega > \Sigma$.
\QED
\vspace{2ex}

The expressions $\mu_k(X)$ given in Proposition~\ref{th:sun-moments}
refer to a SUN variable $X$ with location parameter $\xi=0$.
For the general case with  arbitrary $\xi$, consider the shifted
variable $Y=\xi+X$ and use expressions \ref{eq1}--\ref{eq5} in the appendix.
Another annotation is that Proposition~\ref{th:sun-moments} includes
an expression of $\var{X}$ alternative to \ref{e:sun-var}.

The actual usage of the expressions provided in Proposition~\ref{th:sun-moments}
requires knowledge of the moments of the truncated normal component, $\mu_k(U)$.
In general, these moments are not amenable to explicit treatment, and one must
resort on numerical computations.
As already mentioned, there exists a vast literature concerned with this problem,
and its exhaustive review  would take far too much space.
We therefore indicate only some recent results,
referring the reader to the references quoted therein for earlier developments.
Among the more recent proposals, we  mention the methods for computing
these moments  presented  by \cite{arismendi:2013} and \cite{kan:robo:2017}.
For the latter approach,  there exist publicly available computing routines
written by the authors in the Matlab language.
Of these routines,  a corresponding version is available in the
\texttt{R} computing environment via  either of its packages
\texttt{mnormt} \citep{Rpkg:mnormt} or \texttt{MomTrunc} \citep{Rpkg:MomTrunc}.

We now want to obtain expressions for the Mardia's measures of multivariate
skewness and kurtosis, denoted $\beta_{1,d}$ and $\beta_{2,d}$ in the
original publications of \citet{mardia:1970,mardia:1974},
apart from the symbol $d$ adopted here to denote the dimensionality.
To simplify the algebraic work, it is convenient to work with a suitably
transformed variable, exploiting the invariance properties of Mardia's
measures with respect to nonsingular affine transformations.
For a random variable $Y\sim\SUN_{d,m}(\xi,\Omega,\Delta,\tau,\Gamma)$,
consider again its representation \ref{e:sun-additive2}, and introduce
additionally
\[
  \mu_0 = \Lambda\,\mu_1(U), \quad \mu=\E{Y}=\xi+\mu_0,\quad
  Y = \mu + X_0, \quad X_0=X-\mu_0, \quad U_0=U-\mu_0
\]
where $X$ is as in Proposition~\ref{th:sun-moments}.
Ruling out degenerate cases, $\Sigma=\var{Y}=\var{X_0}$ is non-singular.
Consider then any non-singular $d \times d$ matrix $C$ such that $\Sigma=C\,C\T$;
although not strictly necessary, a common choice is to set $C=\Sigma^{1/2}$,
the unique symmetric positive-definite square root of $\Sigma$.
Next, introduce the standardized variable 
\begin{equation}
  \tilde{Z} = C\inv (Y-\mu) = C\inv \Lambda U_0 + C\inv \Psi^{1/2} V  \sim
  \SUN_{d,m}\left(-C\inv\mu_0, C\inv \Omega (C\inv)\T, C\inv\omega\Delta, \tau, \Gamma\right)
  \label{e:~Z}
\end{equation}
such that
\[  \E{Z}=0, \qquad \var{Z} = I_d  \,.\]

On setting $\tilde\Lambda = C\inv \Lambda$ and $\tilde\Psi=C\inv\Psi(C\inv)\T$,
where $\Lambda$ and $\Psi$ are given in \ref{e:Lambda,Psi}, and a matching definition
of $\tilde\Psi^{1/2}$, we also  note that 
\begin{equation}
   \tilde Z\buildrel d\over = \tilde\Lambda U_0 +  \tilde\Psi^{1/2} V \,.
   \label{e:~Z-equiv}
\end{equation}
where $\buildrel d\over = $ means identically distributed.

The reason for introducing the variable $\tilde Z$ is represented by the following fact.
For a standardized variable $X^*$, say, having zero mean  vector
and identity variance matrix, the Mardia's measures can be conveniently computed
using the expressions given by \cite{koll:sriv:2005}, namely
\[
   \beta_{1,d} = \tr\{\mu_3(X^*)\T \mu_3(X^*)\}
              = \vec\{\mu_3(X^*)\}\T \vec\{\mu_3(X^*)\}\,,
   \qquad\quad
   \beta_{2,d} = \tr\{\mu_4(X^*)\}\,.
\]
The next statement presents the evaluation of these expressions for the SUN 
variable $\tilde{Z}$. 

\begin{proposition}   \label{e:sun-Mardia}
For the random variable $\tilde Z$ specified as in \ref{e:~Z}
or, equivalently, as in \ref{e:~Z-equiv},
the following expected values hold:
\begin{eqnarray*}
\mu_3(\tilde Z) &=&
   (\tilde\Lambda\otimes \tilde\Lambda)\mu_3(U_0)\tilde\Lambda\T,\\
\mu_4(\tilde Z) &=&
   (\tilde\Lambda\otimes \tilde\Lambda)\mu_4(U_0)(\tilde\Lambda\otimes \tilde\Lambda)\T\\
&& + (I_{p^2}+ K_{p})(\tilde\Lambda\Sigma_U \tilde\Lambda\T\otimes\tilde\Psi + \tilde\Psi\otimes \tilde\Lambda\Sigma_U \tilde\Lambda\T + \tilde\Psi\otimes\tilde\Psi)\\
&& + \vec(\tilde\Lambda\Sigma_U \tilde\Lambda\T)\vec(\tilde\Psi)\T + \vec(\tilde\Psi)\vec(\tilde\Lambda\Sigma_U \tilde\Lambda\T)\T  + \vec(\tilde\Psi)\vec(\tilde\Psi)\T,
\end{eqnarray*}
where $\Sigma_U=\mu_2(U_0)=\var{U}$, $\tilde\Lambda = C\inv \Lambda$, 
$\tilde\Psi=C\inv\Psi(C\inv)\T$, and $C\,C\T=\Sigma$.
Moreover, the Mardia's measures of multivariate skewness and kurtosis are
\begin{eqnarray*}
\beta_{1,d}&=&\tr\{(\Lambda\T\Sigma\inv\Lambda\otimes \Lambda\T\Sigma\inv\Lambda)\mu_3(U_0)\Lambda\T\Sigma\inv\Lambda \mu_3(U_0)\T\}\\
&=&\vec\{\mu_3(U_0)\}\T(\tilde\Lambda\T\tilde\Lambda\otimes\tilde\Lambda\T\tilde\Lambda\otimes \tilde\Lambda\T\tilde\Lambda)\vec\{\mu_3(U_0)\},\\
\beta_{2,d}&=&{\rm tr}\{(\Lambda\T\Sigma\inv\Lambda\otimes \Lambda\T\Sigma\inv\Lambda)\mu_4(U_0)\}
 + 2{\rm tr}(\Sigma_U \Lambda\T\Sigma\inv\Lambda){\rm tr}\{\Psi\Sigma\inv\}\\
&& + {\rm tr}\{\Psi\Sigma\inv\}^2
  + 4{\rm tr}\{\Sigma_U\Lambda\T\Sigma\inv\Psi\Sigma\inv\Lambda\}
  + 2\tr\{\Psi\Sigma\inv\Psi\Sigma\inv\},
\end{eqnarray*}
where
\begin{eqnarray*}
\mu_3(U_0)&=&(I_{q^2} + K_q)\{\mu_1(U)\otimes \mu_1(U)\mu_1(U)\T - \mu_1(U)\otimes \mu_2(U)\} - \vec\{\mu_2(U)\}\mu_1(U)\T + \mu_3(U),\\
\mu_4(U_0)&=&-3 \mu_1(U)\mu_1(U)\T\otimes \mu_1(U)\mu_1(U)\T + (I_{q^2} + K_q)\{ \mu_1(U)\mu_1(U)\T\otimes \mu_2(U)\nonumber\\
&& + \mu_2(U)\otimes \mu_1(U)\mu_1(U)\T - (\mu_1(U)\otimes I_q)\mu_3(U)\T  - \mu_3(U)(\mu_1(U)\T\otimes I_q)\}\nonumber\\
&& + \vec\{\mu_2(U)\}(\mu_1(U)\otimes \mu_1(U))\T+ (\mu_1(U)\otimes \mu_1(U))\vec\{\mu_2(U)\}\T + \mu_4(U)\,.
\end{eqnarray*}

\end{proposition}

\noindent{\it Proof.}
The expressions of $\mu_3(\tilde Z)$ and $\mu_4(\tilde Z)$ follow directly from
Proposition~\ref{th:sun-moments}, by using it with the terms $X, \Lambda, \Psi$ 
specified as $\tilde{Z}, \tilde\Lambda, \tilde\Psi$.

Therefore, we concentrate on the derivation of $\beta_{1,p}$ and $\beta_{2,p}$ only.
An algebraically convenient route to obtain these quantities is from the stochastic
representation in \ref{e:~Z}.
Denote by $Y'=\mu+C\,\tilde Z'$ an independent replicate  of
$Y=\mu+C\,\tilde Z$, 
so that the Mardia's measure of skewness can be expressed as
\begin{eqnarray*}
\beta_{1,d}&=&\E{[(Y-\mu)\T\Sigma\inv(Y'-\mu)]^3}\\
           &=&\E{(\tilde Z\T  \tilde Z')^3}.
\end{eqnarray*}
First, introduce the matrices $M_{00}= \Lambda\T\Sigma\inv\Lambda$, 
$M_{01}= \Lambda\T\Sigma\inv\Psi^{1/2}$,
$M_{10}= \Psi^{1/2}\Sigma\inv\Lambda=M_{01}\T$, and
$M_{11}= \Psi^{1/2}\Sigma\inv\Psi^{1/2}$. Then expand
\begin{eqnarray*}
\left(\tilde Z\T \tilde Z'\right)^3
&=& [(U_0\T \Lambda\T (C\inv)^{\top}+V\T \Psi^{1/2} (C\inv)^{\top}) (C\inv \Lambda\,U_0'+C\inv \Psi^{1/2} V')]^3\\
&=& (U_0\T M_{00}U_0'+U_0\T M_{01} V' + V\T M_{01}\T U_0'+V\T M_{11}  V')^3\\
&=& (U_0\T M_{00}U_0'+U_0\T M_{01} V')^3\\
&& + 3(U_0\T M_{00}U_0'+U_0\T M_{01} V')^2(V\T M_{01}\T U_0'+V\T M_{11}  V')\\
&& + 3(U_0\T M_{00}U_0'+U_0\T M_{01} V')(V\T M_{01}\T U_0'+V\T M_{11}  V')^2\\
&& + (V\T M_{01}\T U_0'+V\T M_{11}  V')^3\\
&=& (U_0\T M_{00}U_0')^3 + 3(U_0\T M_{00}U_0')^2(U_0\T M_{01} V')\\
&& + 3(U_0\T M_{00}U_0')(U_0\T M_{01} V')^2 + (U_0\T M_{01} V')^3\\
&& + 3[\{(U_0\T M_{00}U_0')^2+2U_0\T M_{00}U_0'U_0\T M_{01} V'+(U_0\T M_{01} V')^2\} V\T M_{01}\T U_0'\\
&&+\{(U_0\T M_{00}U_0')^2+2U_0\T M_{00}U_0'U_0\T M_{01} V'+(U_0\T M_{01} V')^2\} V\T M_{11}  V']\\
&& + 3[U_0\T M_{00}U_0'\{(V\T M_{01}\T U_0')^2+2V\T M_{01}\T U_0'V\T M_{11}  V'+(V\T M_{11}  V')^2\}\\
&&+U_0\T M_{01} V'\{(V\T M_{01}\T U_0')^2+2(V\T M_{01}\T U_0'V\T M_{11}  V'+(V\T M_{11}  V')^2\}]\\
&& + (V\T M_{01}\T U_0')^3 + 3(V\T M_{01}\T U_0')^2(V\T M_{11}  V')\\
&& + 3(V\T M_{01}\T U_0')(V\T M_{11}  V')^2
+ (V\T M_{11}  V')^3.
\end{eqnarray*}
Take into account that $U_0$ and $U_0'$  are independent and identically distributed (i.i.d.) random
vectors  with mean zero, as well as that $V$ and $V'$ are i.i.d.\ random vectors with spherical normal distribution, so that all the terms involving odd functions of $V$ or $V'$ have zero expectation.
Consider also that $U_0$, $U_0'$, $V$ and $V'$ are mutually independent. 
Then, by taking the expectation and removing the terms with zero mean, we have
\begin{eqnarray*}
 \E{\left(\tilde Z\T  \tilde Z'\right)^3}
&=& \E{\left(U_0\T M_{00}U_0'\right)^3}\\
&=& \E{U_0\T M_{00}U_0'U_0\T M_{00}U_0'U_0\T M_{00}U_0'}\\
&=& \E{(U_0')\T  M_{00}U_0U_0\T M_{00}U_0'(U_0')\T  M_{00}U_0}\\
&=& \E{\tr(M_{00})U_0U_0\T M_{00}U_0'(U_0')\T  M_{00}U_0(U_0')\T }.
\end{eqnarray*}
Now use the equality $\tr(EFGH)=\vec(H\T)\T(G\T\otimes E)\vec(F)$ given
in Lemma~3 of \cite{magn:neud:1986}
with $E=M_{00}U_0U_0\T M_{00}$, $F=U_0'(U_0')\T $,
$G=M_{00}U_0$ and $H=(U_0')\T $, and write
\begin{eqnarray*}
\beta_{1,d}= \E{\left(\tilde Z\T \tilde Z'\right)^3}
&=& \E{\vec(U_0')\T(U_0\T M_{00}\otimes M_{00}U_0U_0\T M_{00})\vec(U_0'(U_0')\T )}\\
&=&\E{ (U_0')\T M_{00}(U_0\T\otimes U_0U_0\T)(M_{00}\otimes M_{00})(U_0'\otimes U_0')}\\
&=&\E{\tr\{(M_{00}\otimes M_{00})(U_0'(U_0')\T \otimes U_0')M_{00}(U_0\T\otimes U_0U_0\T)\}}\\
&=&\tr\left\{(M_{00}\otimes M_{00})\E{U_0'(U_0')\T \otimes U_0'}M_{00}\E{U_0\T\otimes U_0U_0\T)}\right\}\\
&=&\tr\left\{(M_{00}\otimes M_{00})\mu_3(U_0')M_{00} \mu_3(U_0)\T\right\},
\end{eqnarray*}
where $\mu_3(U_0')=\mu_3(U_0)$ since $U_0'$ and $U_0$ are i.i.d.\ variables.

Proceeding in a similar way for the measure of kurtosis, we have
\begin{eqnarray*}
\beta_{2,d}&=&\E{\left[(Y-\mu)\T\Sigma\inv(Y-\mu)\right]^2}\\
&=&\E{(Z\T Z)^2}\\
&=& \E{\left[\left(U_0\T \Lambda\T (C\inv)\T + V\T \Psi^{1/2} (C\inv)\T\right)
    \:\left(C\inv \Lambda\,U_0+C\inv \Psi^{1/2} V\right)\right]^2}\\
&=& \E{(U_0\T M_{00}U_0+U_0\T M_{01}  V + V\T M_{01}\T U_0+V\T M_{11} V)^2}\\
&=& \E{(U_0\T M_{00}U_0+U_0\T M_{01}  V )^2}\\
&&+2\E{(U_0\T M_{00}U_0+U_0\T M_{01}  V)(V\T M_{01}\T U_0+V\T M_{11} V)}\\
&&+\E{(V\T M_{01}\T  U_0+V\T M_{11}  V)^2}\\
&=& \E{(U_0\T M_{00}U_0)^2}+2\E{(U_0\T M_{00}U_0)(U_0\T M_{01} V )}+\E{(U_0\T M_{01}  V )^2}\\
&&+2\E{(U_0\T M_{00}U_0)(V\T M_{01}\T  U_0+V\T M_{11} V)}\\
&&+2\E{(U_0\T M_{01}  V)(V\T M_{01}\T  U_0+V\T M_{11} V)}\\
&&+\E{(V\T M_{01}\T  U_0)^2}+2\E{(V\T M_{01}\T  U_0)(V\T M_{11} V)}+\E{(V\T M_{11} V)^2}\\
&=& \E{(U_0\T M_{00}U_0)^2}+\E{(U_0\T M_{01} V )^2}\\
&&+2\E{(U_0\T M_{00}U_0)(V\T M_{11}V)}+2\E{(U_0\T M_{01}V)(V\T M_{01}\T U_0)}\\
&&+\E{(V\T M_{01}\T  U_0)^2}+\E{(V\T M_{11}V)^2},
\end{eqnarray*}
where the terms with zero expectation have been removed,
namely those associated with odd functions of $V$.
The remaining expected values can be worked out recalling that the powers
of quadratic forms  can be expressed as the trace of matrix products,
combined with properties of the trace of products of matrices,
specifically that
$\tr(K_{pq}(P\T\otimes Q))= \tr(P\T Q)$ as stated by Theorem~3.1, item (xiii),
of \cite{magn:neud:1979as} and, in case $p=q$, $\tr(P\otimes Q)=\tr(P)\tr(Q)$,
$\tr(P\T Q)=\vec(P)\T\vec(Q)$. We then obtain
\begin{eqnarray*}
\beta_{2,d}
&=& \E{\tr(M_{00}U_0U_0\T M_{00}U_0U_0\T)}+2\:\E{\tr(M_{01} VV\T M_{01}\T U_0U_0\T)}\\
  &&+2\:\E{\tr(M_{00}U_0U_0\T)\tr(M_{11}VV\T)}+2\:\E{\tr(M_{01}VV\T M_{01}\T U_0U_0\T)}\\
  &&+\E{\tr(M_{11}VV\T M_{11} VV\T)}\\
&=& \E{\vec(U_0U_0\T)\T(M_{00}\otimes M_{00})\vec(U_0U_0\T)}
  +2\:\E{\vec(U_0U_0\T)\T(M_{01}\otimes M_{01})\vec(VV\T)}\\
&&+2\:\E{\tr(M_{00}U_0U_0\T)\tr(M_{11}VV\T)}
  +2\:\E{\tr(M_{01}  VV\T M_{01}\T U_0U_0\T)}\\
&&+\E{\vec(VV\T)\T(M_{11} \otimes M_{11} )\vec(VV\T)}\\
&=&\tr\left[(M_{00}\otimes M_{00})\E{\vec(U_0U_0\T)\vec(U_0U_0\T)\T}\right]
  +2\:\E{\vec(U_0U_0\T)\T}(M_{01}\otimes M_{01})\E{\vec(VV\T)}\\
&&+2\:\tr\left[M_{00}\E{U_0U_0\T}\right)\tr\left(M_{11}  \E{VV\T}\right]
  +2\:\tr\left[M_{01}  \E{VV\T}M_{01}\T \E{U_0U_0\T}\right]\\
&&+\tr\left[(M_{11}\otimes M_{11})\E{\vec(VV\T)\vec(VV\T)\T}\right]\\
&=&\tr\left\{(M_{00}\otimes M_{00})\mu_4(U_0)\right\}
 +2\vec\{\mu_2(U_0)\}\T(M_{01} \otimes M_{01})\vec\{\mu_2(V)\}\\
&&+2\:\tr\left(M_{00}\mu_2(U_0)\right)\tr\left(M_{11}\mu_2(V)\right)
  +2\:\tr(M_{01} \mu_2(V)M_{01}\T \mu_2(U_0))\\
&&+\tr\left[(M_{11}\otimes M_{11})\mu_4(V)\right] \,.
\end{eqnarray*}

Taking into account Lemma~\ref{th:moments-normal} in an appendix,
we can substitute $\mu_2(V)=I_d$,
$\vec\{\mu_2(V)\}=\vec(I_d)$, $\mu_4(V)=I_{d^2}+K_d+\vec(I_d)\vec(I_d)\T$,
and  $\mu_2(U_0)=\var{U_0}$, leading to
\begin{eqnarray*}
\beta_{2,d}
&=&\tr\left[(M_{00}\otimes M_{00})\mu_4(U_0)\right]
  +2\vec\{\mu_2(U_0)\}\T(M_{01} \otimes M_{01})\vec(I_d)\\
&&+2\:\tr\left(M_{00}\mu_2(U_0)\right)\tr\left(M_{11}\right)
  +2\:\tr(M_{01} M_{01}\T )\\
&&+\tr\left\{(M_{11}\otimes M_{11})(I_{d^2}+K_d)\right\}
  +\tr\left\{(M_{11}\otimes M_{11})\vec(I_d)\vec(I_d)\T\right\}\\
&=&\tr\left\{\mu_4(U_0)(M_{00}\otimes M_{00})\right\}
  +4\:\tr(\mu_2(U_0)M_{01}M_{01}\T )\\
&&+2\:\tr(\mu_2(U_0)M_{00})\tr(M_{11})
  +\tr(M_{11})^2+2\:\tr(M_{11}^2) \,,
\end{eqnarray*}
where $M_{01}M_{01}\T =\Lambda\T \Sigma\inv\Psi\Sigma\inv\Lambda$, 
$\tr(M_{11})=\tr(\Psi\Sigma\inv)$ and $\tr(M_{11}^2)=\tr(\Psi\Sigma\inv\Psi\Sigma\inv)$.
\QED

\section{Other properties} \label{s:other-properties}

\subsection{Log-concavity of the SUN distribution} \label{s:log-concavity}

The SN distribution is known to be log-concave, even in its extended version, ESN;
see \citet{azza:rego:2012aism} for a proof. Since the ESN distribution corresponds
to the SUN with $m=1$, it is natural to investigate the same property for a general
value of $m$.

An often-employed definition of log-concave distribution in the continuous case
requires that the logarithm of its density function is a concave function.
In the more specialized literature, the concept of log-concavity is expressed
via the corresponding probability measure, by requiring that
\begin{equation}
  \pr{\lambda A + (1-\lambda) B} \ge \pr{A}^\lambda\:\pr{B}^{1-\lambda}
  \label{e:log-concave-measure}
\end{equation}
for any two Borel sets $A$ and $B$, and for any $0<\lambda<1$.
For general information on this theme, we refer to Chapter~2  of
\cite{dhar:joag:1988} and  Chapter~4 of \cite{prekopa:1995},
which provide extensive compendia of a vast literature.

Established results ensure  the equivalence of the
definition of log-concavity based on the density function and the one
in \ref{e:log-concave-measure};
see Theorems~4.2.1  of  \cite{prekopa:1995}, and
Theorem~2.8 of \cite{dhar:joag:1988}.
Moreover, also the corresponding distribution function is a log-concave
function; see Theorem~and 4.2.4\,II of \cite{prekopa:1995}.

\begin{proposition}
The SUN distribution is log-concave.
\label{th:sun-log-concavity}
\end{proposition}
\noindent\emph{Proof.}
The proof is based on its additive representation in the form \ref{e:sun-additive2},
which involves the underling variable $Z_0$ indicated in \ref{e:(d+m)-normal'}.
For the multivariate normal distribution, log-concavity is a well-known fact.
Next, recall Theorem~9 of \citet{horrace:2005} which ensures log-concavity
of a normal distribution subject to one-sided truncation.
In our case the truncation  operates on the variable $Z_0=(V\T, W\T)$
in the form $W+\tau>0$. Since $U\equald (W|W+\tau>0)$,
this establishes log-concavity of the distribution of $(V\T, U\T)$.
A variable $Y\sim \SUN_{d,m}(\xi,\Omega,\Delta, \tau, \Gamma)$ can be
obtained from $(V\T, U\T)$ by the affine transformation
\[
   Y = \xi + \pmatrix{\Psi^{1/2} & 0 \cr 0 & \Lambda} \pmatrix{V \cr U}\,.
\]
Preservation of log-concavity after an affine transformation has been
proved by \cite{henn:astr:2006}.
Strictly speaking, their statement refers to a transformation involving
a square matrix,  having dimension $d+m$ in our notation, but it is easy
to see that fact extends to reduced-dimension transformations,
since one can think of a full-rank transformation to an augmented variable
of dimension $d+m$, followed by marginalization to extract the $Y$ component.
Since marginalization preserves log-concavity, as stated for instance
by Theorem~4.2.2 of \cite{prekopa:1995}, this concludes the proof.
\hfill{\scriptsize QED}
\subsection{Conditional density  generated by interval selection}
\label{s:conditional-interval}

For a random variable $Y\sim \SUN_{d,m}(\xi, \Omega, \Delta, \tau, \Gamma)$,
consider a partition of  $Y$ and its associated quantities, as follows
\begin{equation}
   Y = \pmatrix{Y_1 \cr Y_2}, \quad
   \xi= \pmatrix{\xi_1 \cr \xi_2}, \quad
   \Omega= \pmatrix{\Omega_{11} & \Omega_{12} \cr \Omega_{21} & \Omega_{22}}, \quad
   \omega = \pmatrix{\omega_{1} & 0 \cr 0 & \omega_{2}},\quad
   \Delta= \pmatrix{\Delta_1 \cr \Delta_2}
     \label{e:Y-partition}
\end{equation}
where $Y_1$ and $Y_2$ have dimension $d_1$ and $d_2$, with a corresponding partition
for the scaled matrix $\bar\Omega$  which appears in \ref{e:(d+m)-normal}
and \ref{e:sun-pdf}.

In Proposition~2.3.2 of \citet{gonz:etal:2004-inSE}, it is proved
that the conditional distribution of $Y_2$ given that $Y_1=y_1$,
for any vector $y_1\in\Real^{d_1}$, is still of SUN type.
Here, we want to examine another conditional distribution of $Y_2$,
namely the one which arises when the conditioning event on $Y_1$ is
instead an orthant-type interval  of the form $(Y_1+y_1>0)$, where the
inequality sign holds for each variable component, or some similar
orthant-type condition.

\begin{proposition} \label{th:conditional-sun}
If $Y\sim \SUN_{d,m}(\xi, \Omega, \Delta, \tau, \Gamma)$ with elements partitioned
as indicated in \ref{e:Y-partition}, then
\begin{equation}
 (Y_2 | Y_1+y_1>0)\sim \SUN_{d_2,d_1+m}
    \left(\xi_2,\Omega_{22},(\Delta_2,\bar\Omega_{21}),
          \pmatrix{ \tilde{z}_1 \cr \tau},
              \pmatrix{\bar\Omega_{11} &\Delta_1 \cr \Delta_1\T&\Gamma}
    \right)
\label{e:conditional-sun}
\end{equation}
where $\tilde{z}_1=\omega_1\inv(\xi_1+y_1)$ and the inequality sign must be intended
to hold for each component of $Y_1$, if $d_1>1$.
In the case where the inequality sign is reversed, we have
\begin{equation}
 (Y_2 | Y_1+y_1<0)\sim \SUN_{d_2,d_1+m}
    \left(\xi_2,\Omega_{22},(\Delta_2, -\bar\Omega_{21}),
          \pmatrix{ \hat{z}_1 \cr \tau},
              \pmatrix{\bar\Omega_{11} &-\Delta_1\cr -\Delta_1\T&\Gamma}
    \right)
\label{e:conditional-sun-2}
\end{equation}
where $\hat{z}_1=\omega\inv(\xi_1-y_1)$.
\end{proposition}

\noindent\emph{Proof.}
Recall formula \ref{e:sun-affine} for computing the
distribution of an affine transformation of a SUN variable.
Using these transformation rule,  $Z_j=\omega_j\inv(Y_j-\xi_j) \sim
\SUN_{d_j,m}(0, \bar\Omega_{jj}, \Delta_j, \tau, \Gamma)$
 for $j=1,2$.
Then, on setting $z_2=\omega_2\inv(y_2-\xi_2)$, the  conditional
distribution function of $(Y_2|Y_1+y_1>0)$ evaluated at
$y_2\in\Real^{d_2}$  is
\begin{eqnarray}
F_{Y_2}(y_2|Y_1+y_1>0)&=&\pr{Y_2\le y_2\mid Y_1+y_1>0} \nonumber\\
&=&\frac{\pr{Y_1+y_1>0,Y_2\le y_2}}{\pr{ Y_1+y_1>0}}  \nonumber\\
&=&\frac{\pr{\xi_1+\omega_1 Z_1+y_1>0,Y_2\le y_2}}{\pr{\xi_1+\omega_1 Z_1+y_1>0}}  \nonumber\\
&=&\frac{\pr{-Z_1<\tilde{z}_1,Z_2\le z_2}}{\pr{-Z_1<\tilde{z}_1}} \nonumber\\
&=&\frac{F_{-Z_1,Z_2}(\tilde{z}_1,z_2)}{F_{-Z_1}(\tilde{z}_1)}.  \label{e:ratio-cdf}
\end{eqnarray}
where $F_{X}(\cdot)$ denotes the distribution function of a SUN variable $X$,
given by \ref{e:sun-cdf}.
Using again formula \ref{e:sun-affine}, write
\begin{equation}
\pmatrix{-Z_1\cr Z_2}  \sim \SUN_{d_1+d_2,m}
   \left(\pmatrix{0\cr 0},
   \pmatrix{\bar\Omega_{11}&-\bar\Omega_{12}\cr -\bar\Omega_{21}&\bar\Omega_{22}},
   \pmatrix{-\Delta_1\cr \Delta_2},\tau,\Gamma\right)
   \label{e:sun-change-Z1-sign}
\end{equation}
and
\[-Z_1\sim \SUN_{d_1,m}(0,\bar\Omega_{11},-\Delta_1,\tau,\Gamma) \,,\]
so that the two ingredients of \ref{e:ratio-cdf} are
\[F_{-Z_1,Z_2}(\tilde{z}_1,z_2)=\frac{1}{\Phi_m(\tau;\Gamma)}\:
  \Phi_{d_1+d_2+m}\left\{\pmatrix{\tilde{z}_1\cr z_2\cr \tau};
             \pmatrix{\bar\Omega_{11}&-\bar\Omega_{12} & \Delta_1\cr
                      -\bar\Omega_{21}&\bar\Omega_{22} &-\Delta_2\cr
                      \Delta_1\T & -\Delta_2\T & \Gamma}\right\}
\]
 and
\[F_{-Z_1}(\tilde{z}_1)=\frac{1}{\Phi_m(\tau;\Gamma)}\:
  \Phi_{d_1+m}\left\{\pmatrix{ \tilde{z}_1 \cr \tau};
              \pmatrix{\bar\Omega_{11} &\Delta_1\cr\Delta_1\T&\Gamma}\right\}\,,
\]
Taking the ratio of the last two expressions, we obtain
\begin{eqnarray*}
F_{Y_2}(y_2 | Y_1+y_1>0)
&=&\frac{1}{\Phi_{d_1+m}\left\{\pmatrix{\tilde{z}_1 \cr \tau};
 \pmatrix{\bar\Omega_{11} &\Delta_1\cr\Delta_1\T&\Gamma}\right\}}\;
 \Phi_{d_1+d_2+m}\left\{\pmatrix{\tilde{z}_1\cr z_2\cr \tau};
 \pmatrix{\bar\Omega_{11}&-\bar\Omega_{12} & \Delta_1\cr
                      -\bar\Omega_{21}&\bar\Omega_{22} &-\Delta_2\cr
                      \Delta_1\T & -\Delta_2\T & \Gamma}\right\}
\end{eqnarray*}
which is the distribution function of a SUN variable with parameters indicated
in \ref{e:conditional-sun}.

For statement \ref{e:conditional-sun-2}, notice that the event $\{Y_1+y_1<0\}$ coincides
with $\{-Y_1+(-y_1)>0\}$ and apply  \ref{e:conditional-sun} to the distributions of
$(-Y_1,Y_2)$ with $y_1$ replaced by $-y_1$. The distribution of $(-Y_1,Y_2)$ is essentially
given by \ref{e:sun-change-Z1-sign}, up to a change of location and scale.
\hfill{\scriptsize{QED}}
\par\vspace{2ex}
Clearly, the special case of Proposition~\ref{th:conditional-sun} where $m=1$
applies to the extended SN distribution.
By following the same logic of Proposition~\ref{th:conditional-sun}, it is
conceptually simple, although algebraically slightly intricate, to
write the conditional distribution of $(Y_2|Y_1\in I)$ where $I$ is
an event specified by mixed-direction inequalities on the components of $Y_1$.


\appendix
\setcounter{theorem}{0}
    \renewcommand{\thetheorem}{\Alph{section}.\arabic{theorem}}
\setcounter{equation}{0}
    \renewcommand{\theequation}{\Alph{section}.\arabic{equation}}
\section{Appendix}

\subsection{On moments of a variable after a selection}

The following results are presumably well-known.
However, since we are not aware of similarly explicit statements
in the literature,  they are included here.

\begin{lemma} \label{th:mean-by-selection}
For a $m$-dimensional random variable $W$ and an arbitrary Borel
set $A \subseteq \Real^m$,  such that $\pi=\pr{W\in A}$,
consider the associate variables obtained by selection
$U \equald (W|W\in A)$ and $U^c \equald (W|W\not\in A)$.
If the expected value $\E{h(W)}$ exists, for a given function $h$,
then $\E{h(U)}$ and $\E{h(U^c)}$ also exist, and are such that
\begin{equation}
  \E{h(W)}= \E{h(U)}\pi + \E{h(U^c)}(1-\pi) \,.
  \label{e:mean-by-selection}
\end{equation}
\end{lemma}
\noindent\emph{Proof.} Existence of $\E{h(U)}$ and $\E{h(U^c)}$ follows from
the fact that $|h(x)\,I_A(x)|\le |h(x)|$ and $|h(x)\,I_{A^c}(x)|\le |h(x)|$,
where $I_S$ denotes the indicator function of set $S$,
and integrability of $|h(x)|$ is ensured by the existence of $\E{h(W)}$.
The expression in \ref{e:mean-by-selection} follows from law of iterated expectation.
\QED

\begin{proposition} \label{th:var-by-selection}
Under the conditions of Lemma~\ref{th:mean-by-selection}, assume that the variance
matrix $\var{W}$ exists and it is positive semidefinite, written as $\var{W}\ge 0$,
then $\var{U^c}$ and
$\var{U^c}$ also exist, such that
\begin{equation}
\var{W} = \var{U}\pi + \var{U^c}(1-\pi) +
        \left(\E{U}-\E{U^c}\right)\left(\E{U}-\E{U^c}\right)\T\pi(1-\pi) \\
       \label{e:var-by-selection}
\end{equation}
and
\begin{equation}
  \var{W} - \var{U}  \ge  0
  \label{e:var-by-selection-ineq}
\end{equation}
where the inequality sign holds strictly if $\var{W}>0$ and $0<\pi<1$.
\end{proposition}
\noindent\emph{Proof.} Using Lemma~\ref{th:mean-by-selection} with $h$ equal
to the identity function and to the function selecting the generic entry of $W\,W\T$,
write
\begin{eqnarray*}
\E{W} &=& \E{U}\pi + \E{U^c}(1-\pi),\\
\E{WW\T} &=& \E{UU\T}\pi + \E{U^c(U^c)\T}(1-\pi),
\end{eqnarray*}
and then
\begin{eqnarray*}
\var{W} &=& \E{WW\T} - \E{W}\E{W\T}\\
&=& \E{UU\T}\pi + \E{U^c(U^c)\T}(1-\pi)\\
&& - \left[\E{U}\pi + \E{U^c}(1-\pi)\right]\left[\E{U}\pi + \E{U^c}(1-\pi)\right]\T\\
&=& \var{U}\pi + \var{U^c}(1-\pi) + \E{U}\E{U\T}\pi + \E{U^c}\E{(U^c)\T}(1-\pi)\\
&& - \left[\E{U}\pi + \E{U^c}(1-\pi)\right]\left[\E{U}\pi + \E{U^c}(1-\pi)\right]\T\\
&=& \var{U}\pi + \var{U^c}(1-\pi) + \left(\E{U}-\E{U^c}\right)\left(\E{U}-\E{U^c}\right)\T\pi(1-\pi)
\end{eqnarray*}
which proves \ref{e:var-by-selection}. For \ref{e:var-by-selection-ineq}, consider
\begin{eqnarray}
\var{W} - \var{U}
&=& \var{U}(1-\pi) + \var{U^c}(1-\pi) + \left(\E{U}-\E{U^c}\right)\left(\E{U}-\E{U^c}\right)\T\pi(1-\pi)  \nonumber \\
&=&  \left(\var{U} + \var{U^c}\right)(1-\pi) + \left(\E{U}-\E{U^c}\right)\left(\E{U}-\E{U^c}\right)\T\pi(1-\pi)  \label{e:varW-varU}\\
&\ge & 0 \nonumber
\end{eqnarray}
since the two summands of \ref{e:varW-varU} are non-negative definite matrices.\par
Consider now the case when $\var{W}>0$ and $0<\pi<1$.
For an arbitrary vector $a\in\Real^m$,  define $W_a=a\T W$ and $U_a= a\T U$.
Provided $a\not=0$, the condition  $\var{W}>0$ ensures that $\var{W_a}= a\T\var{W}a>0$,
which means that $W_a$ is a non-degenerate variable.
To show that also $U_a$ is a non-degenerate variable, assume that the opposite holds,
which means that there exists a vector $a$ such that $U_a=a\T U \equiv b$,
for some constant $b$. Then
 \begin{eqnarray*}
1=\pr{a\T U=b}&=&\pr{a\T W=b | W\in A}\\
&=&\frac{\pr{a\T W=b, W\in A}}{\pr{W\in A}}\leq \frac{\pr{a\T W=b}}{\pr{W\in
A}}=\frac{0}{\pi}=0
\end{eqnarray*}
where the last equality uses the condition $\pi>0$.
Since we have obtained a contradiction, then $U_a$ cannot be degenerate
and $\var{U}>0$. By a similar argument and the condition $\pi<1$, we can
establish that $\var{U^c}>0$. Therefore the term $\var{U} + \var{U^c}$ in
\ref{e:varW-varU} is positive definite, while the final summand of
\ref{e:varW-varU} is at least positive semidefinite, implying
that $\var{W} - \var{U} >0$.
\QED

\subsection{On moments of multivariate normal variables}
\begin{lemma}  \label{th:moments-normal}
If $V_0\sim \N_r(0,I_r)$, then
\[
\begin{array}{rrcl}
(i)  &\E{V_0\otimes V_0\T} &=&\E{V_0V_0\T}=I_r,\\
(ii) &\E{V_0V_0\T\otimes V_0} &=&\E{V_0\otimes V_0V_0\T}=0,\\
(iii)&\E{V_0V_0\T\otimes V_0V_0\T} &=& I_{r^2}+K_r+\vec(I_r)\vec(I_r)\T,\\
(iv) &\var{\vec(V_0V_0\T)}
     &=& I_{r^2}+K_r.
\end{array}
\]
\end{lemma}
The proof of statements (i) and (ii) is direct.
For (iii) and (iv), see for instance Theorem~4.1\,(i) and Lemma 4.1\,(ii)
of \cite{magn:neud:1979as}.


\subsection{On moments of the sum of two independent multivariate
random variables}

\begin{proposition}  \label{th:moments-of-sum}
Let $X=AU+BV$, where $A\in\mathbb{R}^{p\times q}$ and $B\in\mathbb{R}^{p\times r}$ are constant
matrices, and $U\in\Real^q$ and $V\in\Real^r$ are independent random vector.
If the required moments exist, then
\begin{eqnarray*}
\mu_1(X)&=&A\mu_1(U)+B\mu_1(V)=A\E{U}+B\E{V},\\
\mu_2(X)&=&A\mu_2(U)A\T+A\mu_1(U)\mu_1(V)\T B\T+B\mu_1(V)\mu_1(U)\T A\T + B\mu_2(V)B\T,\\
\mu_3(X)
&=&(A\otimes A)\mu_3(U)A\T+(A\otimes A)\vec\{\mu_2(U)\}\mu_1(V)\T B\T\\
&&+(I_{p^2}+K_p)(A\otimes B)\left\{(\mu_2(U)\otimes \mu_1(V))A\T
                                  +(\mu_1(U)\otimes \mu_2(V))B\T\right\}\\
&&+(B\otimes B)\vec\{\mu_2(V)\}\mu_1(U)\T A\T + (B\otimes B)\mu_3(V)B\T,\\
\mu_4(X)&=&(A\otimes A)\mu_4(U)(A\otimes A)\T\\
&&+(A\otimes A)\mu_3(U)(I_q\otimes \mu_1(V))\T(A\otimes B)\T\\
&&+K_{p}(A\otimes A)\mu_3(U)(I_q\otimes \mu_1(V))\T(A\otimes B)\T K_{p}\\
&&+(A\otimes B)(I_q\otimes \mu_1(V))\mu_3(U)\T(A\otimes A)\T\\
&&+K_{p}(A\otimes B)(I_q\otimes \mu_1(V))\mu_3(U)\T(A\otimes A)\T K_{p}\\
&&+(A\otimes B)(\mu_2(U)\otimes \mu_2(V))(A\otimes B)\T\\
&&+(A\otimes B)(\mu_2(U)\otimes \mu_2(V))(A\otimes B)\T K_{p}\\
&&+K_{p}(A\otimes B)(\mu_2(U)\otimes \mu_2(V))(A\otimes B)\T\\
&&+K_{p}(A\otimes B)(\mu_2(U)\otimes \mu_2(V))(A\otimes B)\T K_{p}\\
&&+(A\otimes A)\vec\{\mu_2(U)\}\vec\{\mu_2(V)\}\T(B\otimes B)\T\\
&&+(B\otimes B)\vec\{\mu_2(V)\}\vec\{\mu_2(U)\}\T(A\otimes A)\T\\
&&+(A\otimes B)(\mu_1(U)\otimes I_r)\mu_3(V)\T(B\otimes B)\T\\
&&+K_{p}(A\otimes B)(\mu_1(U)\otimes I_r)\mu_3(V)\T(B\otimes B)\T K_{p}\\
&&+(B\otimes B)\mu_3(V)(\mu_1(U)\otimes I_r)\T(A\otimes B)\T\\
&&+K_{p}(B\otimes B)\mu_3(V)(\mu_1(U)\otimes I_r)\T(A\otimes B)\T K_{p}\\
&&+(B\otimes B)\mu_4(V)(B\otimes B)\T.
\end{eqnarray*}
\end{proposition}
\noindent{\it Proof:}
The proof of $\mu_1(X)$ is trivial. To obtain the other moments, first note that
\begin{eqnarray*}
XX\T&=&(AU+BV)(U\T A\T + V\T B\T)\\
&=&AUU\T A\T+AUV\T B\T+BVU\T A\T + BVV\T B\T,
\end{eqnarray*}
leading to $\mu_2(X)$. Also,
\begin{eqnarray*}
X\otimes XX\T&=&(AU+BV)\otimes(AUU\T A\T+AUV\T B\T+BVU\T A\T + BVV\T B\T)\\
&=&(AU\otimes AUU\T A\T)+(AU\otimes AUV\T B\T)\\
&&+(AU\otimes BVU\T A\T) + (AU\otimes BVV\T B\T)\\
&&+(BV\otimes AUU\T A\T)+(BV\otimes AUV\T B\T)\\
&&+(BV\otimes BVU\T A\T) + (BV\otimes BVV\T B\T)\\
&=&(A\otimes A)(U\otimes UU\T)A\T+(A\otimes A)(U\otimes UV\T)B\T\\
&&+(A\otimes B)(U\otimes VU\T)A\T + (A\otimes B)(U\otimes VV\T)B\T\\
&&+(B\otimes A)(V\otimes UU\T)A\T+(B\otimes A)(V\otimes UV\T)B\T\\
&&+(B\otimes B)(V\otimes VU\T)A\T + (B\otimes B)(V\otimes VV\T)B\T\\
&=&(A\otimes A)(U\otimes UU\T)A\T+(A\otimes A)(U\otimes U)V\T B\T\\
&&+(A\otimes B)(UU\T\otimes V)A\T + (A\otimes B)(U\otimes VV\T)B\T\\
&&+(B\otimes A)(V\otimes UU\T)A\T+(B\otimes A)(VV\T\otimes U) B\T\\
&&+(B\otimes B)(V\otimes V)U\T A\T + (B\otimes B)(V\otimes VV\T)B\T\\
&=&(A\otimes A)(U\otimes UU\T)A\T+(A\otimes A)(U\otimes U)V\T B\T\\
&&+(I_{p^2}+K_p)(A\otimes B)(UU\T\otimes V)A\T + (I_{p^2}+K_p)(A\otimes B)(U\otimes VV\T)B\T\\
&&+(B\otimes B)(V\otimes V)U\T A\T + (B\otimes B)(V\otimes VV\T)B\T,
\end{eqnarray*}
where we have used
\begin{eqnarray*}
(B\otimes A)(V\otimes UU\T)&=&(B\otimes A)K_{rq}(UU\T\otimes V)=K_p(A\otimes B)(UU\T\otimes V);\\
(B\otimes A)(VV\T\otimes U)&=&(B\otimes A)K_{rq}(V\otimes UU\T)=K_p(A\otimes B)(U\otimes VV\T).
\end{eqnarray*}
This leads to $\mu_3(X)$. Finally,
\begin{eqnarray*}
XX\T\otimes XX\T&=&(AUU\T A\T+AU V\T B\T+BVU\T A\T+BVV\T B\T)\\
&&\otimes(AUU\T A\T+AUV\T B\T+BVU\T A\T+BVV\T B\T)\\
&=&(AUU\T A\T\otimes AUU\T A\T)+(AUU\T A\T\otimes AUV\T B\T)\\
&&+(AUU\T A\T\otimes BVU\T A\T)+(AU U\T A\T\otimes BVV\T B\T)\\
&&+(AUV\T B\T\otimes AUU\T A\T)+(AUV\T B\T\otimes AUV\T B\T)\\
&&+(AUV\T B\T\otimes BVU\T A\T)+(AUV\T B\T\otimes BVV\T B\T)\\
&&+(BVU\T A\T\otimes AUU\T A\T)+(BVU\T A\T\otimes AUV\T B\T)\\
&&+(BVU\T A\T\otimes BVU\T A\T)+(BVU\T A\T\otimes BVV\T B\T)\\
&&+(BVV\T B\T\otimes AUU\T A\T)+(BVV\T B\T\otimes AUV\T B\T)\\
&&+(BVV\T B\T\otimes BVU\T A\T)+(BVV\T B\T\otimes BVV\T B\T)\\
&=&(AUU\T A\T\otimes AUU\T A\T)+(AUU\T A\T\otimes AUV\T B\T)\\
&&+(AUU\T A\T\otimes BVU\T A\T)+(AU U\T A\T\otimes BVV\T B\T)\\
&&+K_p(AUU\T A\T\otimes AUV\T B\T)K_p+(AUV\T B\T\otimes AUV\T B\T)\\
&&+(AUU\T A\T\otimes BVV\T B\T)K_p+(AUV\T B\T\otimes BVV\T B\T)\\
&&+K_p(AUU\T A\T\otimes BVU\T A\T)K_p+(AUU\T A\T\otimes BVV\T B\T)K_p\\
&&+(BVU\T A\T\otimes BVU\T A\T)+(BVU\T A\T\otimes BVV\T B\T)\\
&&+K_p(AUU\T A\T\otimes BVV\T B\T)K_p+K_p(AUV\T B\T\otimes BVV\T B\T)K_p\\
&&+K_p(BVU\T A\T\otimes BVV\T B\T)K_p+(BVV\T B\T\otimes BVV\T B\T),
\end{eqnarray*}
where we used the fact that, if $C$ and $D$ are $p\times q$ and $s\times t$ matrices, respectively,
then the following equalities hold: $(D\otimes C)=K_{sp}(C\otimes D)K_{qt}$,
$K_{ps}(D\otimes C)=(C\otimes D)K_{qt}$, $(D\otimes C)K_{tq}=K_{sp}(C\otimes D)$.

After rearranging common terms, we obtain
\begin{eqnarray*}
XX\T\otimes XX\T
&=&(A\otimes A)(UU\T\otimes UU\T)(A\T\otimes A\T)\\
&&+(A\otimes A)(UU\T\otimes U)(I_q\otimes V\T)(A\T\otimes B\T)\\
&&+K_{p}(A\otimes A)(UU\T\otimes U)(I_q\otimes V\T)(A\T\otimes B\T)K_{p}\\
&&+(A\otimes B)(I_q\otimes V)(UU\T\otimes U\T)(A\T\otimes A\T)\\
&&+K_{p}(A\otimes B)(I_q\otimes V)(UU\T\otimes U\T)(A\T\otimes A\T)K_{p}\\
&&+(A\otimes B)(UU\T\otimes V V\T)(A\T\otimes B\T)\\
&&+(A\otimes B)(UU\T\otimes VV\T)(A\T\otimes B\T)K_{p}\\
&&+K_{p}(A\otimes B)(UU\T\otimes VV\T)(A\T\otimes B\T)\\
&&+K_{p}(A\otimes B)(UU\T\otimes VV\T)(A\T\otimes B\T)K_{p}\\
&&+(A\otimes A)(U\otimes U)(V\T\otimes V\T)(B\T\otimes B\T)\\
&&+(B\otimes B)(V\otimes V)(U\T\otimes U\T)(A\T\otimes A\T)\\
&&+(A\otimes B\T)(U\otimes I_r)(V\T\otimes VV\T)(B\T\otimes B\T)\\
&&+K_{p}(A\otimes B)(U\otimes I_r)(V\T\otimes VV\T)(B\T\otimes B\T)K_{p}\\
&&+(B\otimes B)(V\otimes V V\T)(U\T\otimes I_r)(A\T\otimes B\T)\\
&&+K_{p}(B\otimes B)(V\otimes VV\T)(U\T\otimes I_r)(A\T\otimes B\T)K_{p}\\
&&+(B\otimes B)(VV\T\otimes VV\T)(B\T\otimes B\T)
\end{eqnarray*}
whose expectation is $\mu_4(X)$. \QED

\vspace{2ex}\par
The above results simplify considerably if one of the two variables is
symmetric about the origin, hence with  null odd-order moments.
This the case of interest for us, since one of summands in each of
\ref{e:sun-additive} and \ref{e:sun-additive2} is of this type.
The next statement is the pertaining corollary of
Proposition~\ref{th:moments-of-sum}.

\begin{proposition} \label{th:moments-of-sum0}
Let $X=AU+BV$, where $A\in\mathbb{R}^{p\times q}$ and $B\in\mathbb{R}^{p\times r}$ are constant matrices, and $U\in\Real^q$ and $V\in\Real^r$ are independent random vectors.
If $\mu_1(V)$ and $\mu_3(V)$ are zero, then
\begin{eqnarray*}
\mu_1(X)&=&A\mu_1(U),\\
\mu_2(X)&=&A\mu_2(U)A\T + B\mu_2(V)B\T,\\
\mu_3(X)&=&(A\otimes A)\mu_3(U)A\T + (I_{p^2} + K_p)(A\otimes B)(\mu_1(U)\otimes \mu_2(V))B\T\\
&& + (B\otimes B)\vec\{\mu_2(V)\}\mu_1(U)\T A\T,\\
\mu_4(X)&=&(A\otimes A)\mu_4(U)(A\otimes A)\T + (I_{p^2}+K_p)\{(A\otimes B)(\mu_2(U)\otimes \mu_2(V))(A\otimes B)\T\\
&& + (B\otimes A)(\mu_2(V)\otimes \mu_2(U))(B\otimes A)\T\} + (A\otimes A)\vec\{\mu_2(U)\}\vec\{\mu_2(V)\}\T(B\otimes B)\T\\
&& + (B\otimes B)\vec\{\mu_2(V)\}\vec\{\mu_2(U)\}\T(A\otimes A)\T  + (B\otimes B)\mu_4(V)(B\otimes B)\T
\end{eqnarray*}
provided the required moments exist. Moreover,
\begin{eqnarray*}
\var{X}&=&\mu_2(X)-\mu_1(X)\mu_1(X)\T\\
&=&A\,\var{U}A\T + B\,\var{V}B\T,\\
  \cov{X\otimes X,X)}&=&\mu_3(X)-\vec\{\mu_2(X)\}\mu_1(X)\T\\
&=&(A\otimes A)\cov{U\otimes U,U})A\T + (I_{p^2} + K_p)(A\otimes B)(\mu_1(U)\otimes \mu_2(V))B\T,\\
\var{X\otimes X}&=&\mu_4(X)-\vec\{\mu_2(X)\}\vec\{\mu_2(X)\}\T\\
&=&(A\otimes A)\var{U\otimes U}(A\otimes A)\T + (B\otimes B)\var{V\otimes V}(B\otimes B)\T\\
&& + (I_{p^2}+K_p)\big\{(A\otimes B)\left(\mu_2(U)\otimes \mu_2(V)\right)\,(A\otimes B)\T\\
&& + (B\otimes A)\,\left(\mu_2(V)\otimes \mu_2(U)\right)\,(B\otimes A)\T\big\}.
\end{eqnarray*}
\end{proposition}


For the shifted variable $Y=\xi+ X$, where $\xi$ is an arbitrary $p$-vector, the moments are:
\begin{eqnarray}
\mu_1(Y)&=&\xi + \mu_1(X),\label{eq1}\\
\mu_2(Y)&=&\xi\xi\T + \xi \mu_1(X)\T + \mu_1(X)\xi\T + \mu_1(X)\mu_1(X)\T,\label{eq2}\\
\mu_3(Y)&=&\xi\otimes \xi\xi\T + (\xi\otimes \xi)\mu_1(X)\T
+ (I_{p^2} + K_p)\{\xi\xi\T\otimes \mu_1(X) + \xi\otimes \mu_2(X)\}\nonumber\\
&& + m_2(X)\xi\T + \mu_3(X),\label{eq4}\\
\mu_4(Y)&=&\xi\xi\otimes\xi\xi\T + (I_{p^2} + K_p)\{\xi\xi\T\otimes\xi \mu_1(X)\T + \xi\xi\T\otimes \mu_1(X)\xi\T\nonumber\\
&& + \xi\xi\T\otimes \mu_2(X)+ \mu_2(X)\otimes\xi\xi\T + (\xi\otimes I_p)\mu_3(X)\T  + \mu_3(X)(\xi\T\otimes I_p)\}\nonumber\\
&&+ \mu_2(X)(\xi\otimes\xi)\T+ (\xi\otimes\xi)\mu_2(X)\T + \mu_4(X).\label{eq5}
\end{eqnarray}



\end{document}